\documentclass[a4paper,10pt]{article}
\usepackage{anysize}
\marginsize{3.0cm}{3.0cm}{2.5cm}{1.5cm}
\usepackage{amsmath,amssymb,amsthm,ulem}
\usepackage{latexsym,amsfonts,amsmath}
\usepackage{graphicx,psfrag}
\usepackage{epstopdf}
\usepackage{pstricks}
\usepackage{amssymb}
\usepackage{amsthm}
\usepackage{float}

\def\RR{\mathbb{R}}
\def\QQ{\mathbb{Q}}
\def\II{\mathbb{I}}

\newtheorem{condition}{Condition}[section]{\bfseries}{\itshape}
\newtheorem{remark}{Remark}[section]{\bfseries}{\itshape}
\newtheorem{lemma}{Lemma}[section]{\bfseries}{\itshape}

\begin{document}

\author{
Alexey Piunovskiy \\ Bakhti Vasiev \\ \  \\
Department of Mathematical Sciences, University of Liverpool, L69 7ZL, UK.\\ \texttt{piunov@liverpool.ac.uk, bnvasiev@liverpool.ac.uk}
}
\title{Modelling Ethnogenesis}
\date{}

\maketitle

\maketitle

\begin{abstract} Following the ideas of L.N.Gumilev, we introduce the mathematical model of ethnogenesis which describes the dynamics of subgroups in the developing polity in terms of ordinary differential equations. The bust dynamics associated with the rise and fall of civilisations is modelled as an excitation process, which is the non-linear phenomenon, well known in mathematical biology. We consider deterministic as well as the stochastic version of the model. We also expand the model to study the interaction between two polities undergoing ethnogenesis. Investigation is performed using analytical methods as well as numerical integration (i.e. MATLAB simulation).
\end{abstract}
\begin{tabbing}
\small \hspace*{\parindent}  \= {\bf Keywords:} Ethnogenesis, Population Dynamics, Dynamical System, Excitation \\
\> {\bf AMS 2020 subject classification:} \= 37N25, 37M05, 92B05, 92D15, 92D25
\end{tabbing}

\section{Introduction}

History of mankind in the last 5 thousand years can be viewed in terms of rises and falls of civilizations. L.N.Gumilev considered each civilization as a manifestation of particular ethnos, which, under certain conditions, appears and builds civilization, but in course of time gets old and dies, causing for the associated civilization to disappear. According to  Gumilev's theory \cite{Gumilev2,Gumilev1} the driving force for formation of  a new ethnos comes from a certain group of people, whom he calls "passionaries" (the term which is also used in \cite{Turchin1}) or people with drive. This group is considered as a fraction of population   habitating   certain geographic territory, who express high level of passion and lead their fellows (i.e. tribesmen) forming the rest of the population to expand and to build new society (civilization).

\begin{figure}[ht]
	\centering
	\includegraphics[width=0.7\textwidth]{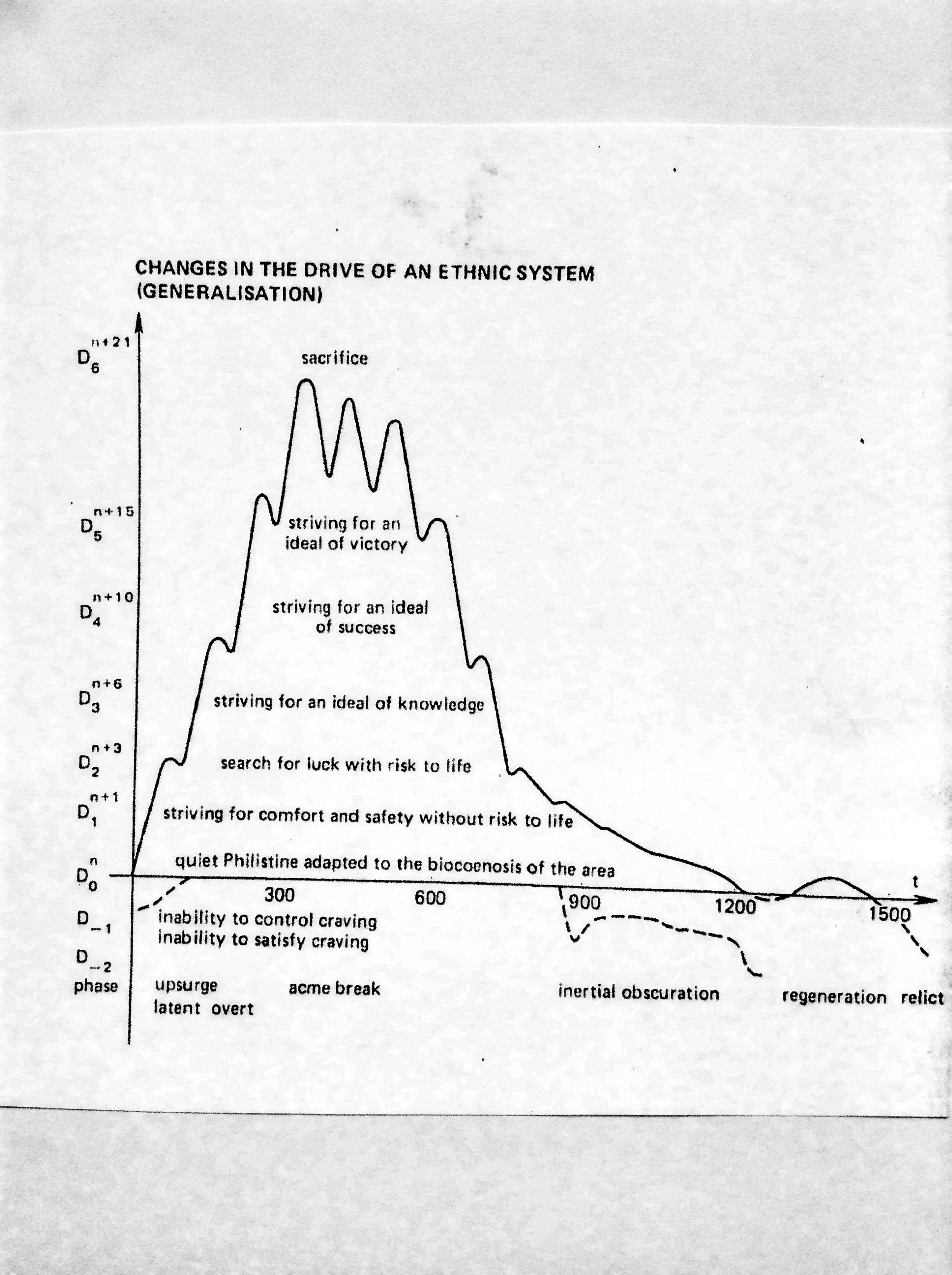}
	\caption{Graphical illustration of the evolution of an ethnos (civilization) \cite[p.240]{Gumilev2}.}
\label{fig:Gum}
\end{figure}

According to Gumilev \cite{Gumilev2,Gumilev1} the formation and death of civilization can be described by the dynamics of civilisation's "passionary tension"  or drive, which he illustrated by a "bust" curve shown in Figure \ref{fig:Gum}.  It starts with a growing phase (rise of civilization) followed by  plateau  (ackmatic phase), fast decline (breaking phase), slow decline (inertial obscuration) and low level tail (obscuration  or regeneration-relict). The entire process, which according to Gumilev takes about 15 centuries, can be considered as a response to a disturbance, caused by initiation of a small fraction of passionaries. Such response is known in physiology as excitation \cite[pp.239-242]{MurrayI}, that is when a small perturbation to the system results in a full-sized response. 

The aim of this study is to develop a mathematical model which would explain the bust dynamics (which is evident from  Figure \ref{fig:Gum}) exhibited in course of ethnogenesis. We note that mathematical study of ethnogenesis was also aimed by other researchers.  However, the mathematical models developed in their works (see, for example, \cite{Guz,Turchin1}, do not appropriately reproduce the dynamics of ethnogenesis, as they didn't consider the association of the bust dynamics with the type of system's nonlinearity known as excitability. The excitable dynamics in the model describing the ethnogenesis is the main theme of the current work. We will build a mathematical model of ethnogenesis, which is based on the statements underlying  Gumilev's theory.  
\begin{itemize}
\item According to Gumilev there are three main subgroups in the population having different behavioural patterns and affecting the evolution of the ethnos. The driving force for  the growth of populations and further improvements resulting in the rise of a civilization comes from passionaries, whose idealistic motivations are grounded on altruism. The bulk of population is represented by harmonious  individuals  who  work on preservation of the current state of the ethnos. There is also a destructive group of individuals, called subpassionaries (or people with negative drive -- vargants, soldier tramps, degenerates), who are as active as passionaries but whose actions are based on egoism rather than altruism. 

\item The formation of a new ethnos (or civilisation) is associated with initiation of passionaries. Gumilev describes conditions under which passionaries appear and take over the population, but we will not go into details of these  conditions and will postulate that at a certain time, a small fraction of population is already represented by passionaries.
\end{itemize}

We will design a few models describing the ethnogenesis. In the first model we will consider the population as consisting of two subgroups, namely, the passionaries and the remainder of the population. This two-variable model will let us to identify  possible interactions between these subgroups which allow the bust dynamics in course of ethnogenesis. In the second version of the model we will consider all three subgroups and analyse the interactions between them which is consistent with the observed dynamics during ethnogenesis. In the follow up steps of our research we will use the three-variable model for the study of the impact of noisy environment to the ethnogenesis. Particularly we show that the noise is amplified by nonlinearities in the model. Finally, we extend  the three variable model (with noise) to consider interaction between two ethnogenetic processes taking place simultaneously, with a certain time lag. This allows us to model conflicting civilisations and to identify conditions when one of them  takes  over the other.

\section{Two-Variable Model}

As a starting point we will consider the two-variable model: 
\begin{equation}\label{eq.2var}
\left\{
\begin{array}{l}
\dot{x}=xf(x,y) \\
\dot{y}=yg(x,y)
\end{array}
\right.
\end{equation}
where variable $x$ represents the size of  the subpopulation formed by passionaries while $y$ is the size of  the remaining population (which includes both harmonious people and subpassionaries). Equation \eqref{eq.2var} is commonly used for modelling population dynamics in biology, where $x$ and $y$ are considered as the sizes of two biological species. Using linear approximation of functions $f(x,y)$ and $g(x,y)$ we get: 
\begin{equation}\label{eq.linear}
\left\{
\begin{array}{l}
\dot{x}=x(a_0+a_1x+c_1y); \\
\dot{y}=y(b_0+b_1y+d_1x),
\end{array}
\right.
\end{equation}
which is a generalised representation of   the Lotka-Volterra model \cite{Lotka1}. Commonly this model is considered under the following condition for model parameters: $a_0>0$, $b_0>0$ while $a_1<0$ and $b_1<0$ preventing unlimited growth of populations. After nondimensionalisation this equation is commonly transferred into 
\begin{equation}\label{eq.nondim}
\left\{
\begin{array}{l}
\dot{x}=x(1-x+\beta_1y); \\
\dot{y}=\gamma y(1-y+\beta_2x),
\end{array}
\right.
\end{equation}
where $\gamma$ defines the relative rate of change of $y$ with respect to $x$ \cite[p.119]{MurrayI}. If $\beta_{1,2}=0$ the two equations are detached and the both species exhibit  the  logistic growth. 
Furthermore, depending on the signs of these two parameters the model reproduces three types of interactions between populations $x$ and $y$, namely, predator-pray ($\beta_1>0$, $\beta_2<0$), symbioses ($\beta_1>0$, $\beta_2>0$) and competition ($\beta_1<0$, $\beta_2<0$) \cite{Lotka1}. The model allows four equilibria and in case of competitive Lotka-Volterra model ($\beta_1<0$, $\beta_2<0$) all four equilibria are meaningful and correspond to non-negative sizes of populations.  These are  the trivial equilibrium $(x=0,~y=0)$, the extinction of $y$-population ($x\ne0$, $y=0$), the extinction of $x$-population ($x=0$, $y\ne0$) and the co-existence ($x\ne0$, $y\ne0$). Depending on model parameters, the solution of the system converges either to the co-existence of species or to the case when one of them becomes extinct \cite[pp.104-126]{MurrayI}. 

In our case we would like for one of the variables, say $x$ to represent the size of  the subpopulation formed by passionaries, while $y$ is the size of the remaining part of population (which would include harmonious people and subpassionaries). As we expect to observe the excitation dynamics, and the excitation in biology is known to be a non-linear process, we will need more than just a linear expansion of functions $f(x,y)$ and $g(x,y)$ in \eqref{eq.2var}. The simplest way is to add one quadratic term into the first equation which transforms our system to  
\begin{equation}\label{eq.nonlinear}
\left\{
\begin{array}{l}
\dot{x}=x(a_0+a_1x+a_2x^2+c_1y); \\
\dot{y}=y(b_0+b_1y+d_1x).
\end{array}
\right.
\end{equation}
There are up to six equilibria in this system of which two are always real:   
$$ (x_1, y_1)=(0,0),  (x_2, y_2)=(0,-b_0/b_1). $$ 
These two equilibria are located on the vertical axis. Two more equilibria (if real) are located on the horizontal axis: 
\begin{equation}\label{eq.roots}
(x_{3,4}, y_{3,4})=\left(\frac{-a_1\pm\sqrt{a_1^2-4a_0a_2}}{2a_2}, 0\right).
\end{equation}
And, finally, two remaining equilibria (if real) are represented by the points of intersection of the line $b_0+b_1y+d_1x=0$ with parabola $a_0+a_1x+a_2x^2+c_1y=0$. To set an excitable kinetics in the system \eqref{eq.nonlinear} we make sure that the equilibrium $(x_2, y_2)$ is stable and located in the vicinity of the parabola $a_0+a_1x+a_2x^2+c_1y=0$. The excitable dynamics becomes more evident  after nondimensionalisation of the system \eqref{eq.nonlinear} (so that $a_2=-1$) and transferring it into the following form: 
\begin{equation}\label{eq.nonlinear1}
\left\{
\begin{array}{l}
\dot{x}=x((1-x)(x-\alpha)+\beta_1 (y-y_0)); \\
\dot{y}=\gamma y(y_0- y+\beta_2 x),
\end{array}
\right.
\end{equation}
where parameters $\gamma$, $\beta_1$ and $\beta_2$ have the same meaning as those in the equation \eqref{eq.nondim}, new parameter $\alpha$ defines the excitation threshold of the system and $y_0$ defines the location of the equilibrium $(x_2,y_2)=(0,y_0)$ which is stable if $\beta_1$ and $\beta_2$ have opposite signs. Note, that the parameters $\alpha$, $y_0$ and $\gamma$ should all be positive. It looks that, in  any case, the system (\ref{eq.nonlinear1}) with positive initial conditions admits a unique bounded solution on the infinite horizon $[0,\infty)$. For $\beta_1\le 0$, $\beta_2\ge 0$, this  will follow from the investigation of the three-variable model (\ref{eq.three-var1}). 
Concerning the  equilibria in the system \eqref{eq.nonlinear1}, we note  that the trivial steady state $(x_1,y_1)=(0,0)$ is unstable (saddle). Furthermore, if the parabola $(1-x)(x-\alpha)+\beta_1 (y-y_0)$ has real and non-negative roots, then the system \eqref{eq.nonlinear1} has two meaningful equilibria located on the horizontal axis (same as given by \eqref{eq.roots}) and for the concave up parabola the equilibrium which is closer to the origin (smaller $x$-coordinate) is the unstable node, while the other one is a saddle. For simplicity, we consider the cases when the nullcline represented by the parabola $(1-x)(x-\alpha)+\beta_1 (y-y_0)=0$ doesn't intersect the one given by the line $y_0- y+\beta_2 x=0$ and therefore we don't have any extra equilibria. 

\begin{figure}[ht]
\centering
\includegraphics[width=1.\textwidth]{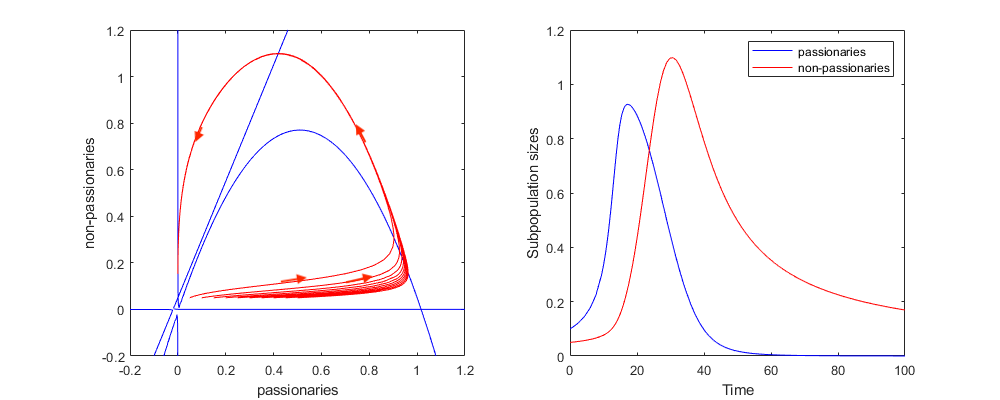}
\caption{Illustration of the excitable dynamics observed in the system \eqref{eq.nonlinear1}. Panel A: nullclines $\dot x=0$ and $\dot y=0$ are shown in blue and a set of phase trajectories with starting points  $x=0.05,0.1, ..., 0.5$ and $y=y_0$ are shown in red. Panel B: time evolution of both variables, $x$ (blue line) and $y$ (red line) from the initial condition  $x=0.1$, $y=y_0$. Parameters values: $\alpha=0.02$, $y_0=0.05$, $\beta_1=-1/3$, $\beta_2=2.5$, $\gamma=0.1$. }
\label{fig:Null}
\end{figure} 

The dynamics in the system \eqref{eq.nonlinear1} is  illustrated by Figure \ref{fig:Null}. Null-clines of the system are shown in blue on panel A. These null-clines indicate the excitable nature of the system \eqref{eq.nonlinear1} and this is illustrated by a set of phase trajectories (shown in red) coursed by the perturbation of the system from its stable equilibrium  $(x_2,y_2)=(0,y_0)$. Any perturbation from this state results  in the relaxation of the system back to this equilibrium. However if the perturbation is above certain threshold, for example,   if  the initial value of $y$ is equal to $y_0$ and the initial value of $x$ is above $\alpha$,
 then the perturbation increases further before the system relaxes back to the stable equilibrium  $(x_2,y_2)$.

The dynamics of variables $x$ and $y$ over time for one of the phase trajectories (starting from the point  $x_0=0.1$ and $y_0$) is shown in panel B. Here we see that  the both variables increase and then decrease over time. We note a relatively fast dynamics of variable $x$ (passionaries) with the duration of the spike being about 60 time unites. We also note a slow relaxation of variable $y$ (the rest of population) which gets back to its equilibrium value with relaxation time of about 100 time units. 

In order to scale the model time units to real time we take into account that according to Gumilev \cite{Gumilev2} the duration of passionary spike is about 900 years which should correspond to 60 time units in the model. Thus one model time unit corresponds to 15 years. Furthermore, according to Gumilev \cite{Gumilev2} the fraction of passionaries can only be up to 5-7\% of the entire population. For the set of parameters values used to produce Fig. \ref{fig:Null} both variables have values roughly in the range (0,1). If we consider variable  $x$ as representing 1\% of  $y$,  the total number of  nonpassionaries, then at the top of the spike (which takes place at  $t\approx15$ and where $x\approx 1$ and $y\approx 0.2$) the passionaries constitute about 5\% of the entire population.
In other words, $x=1$ corresponds, e.g., to the absolute value of $K=10,000$ passionaries, while $y=1$ corresponds to $100K=1,000,000$ nonpassionaries. One can certainly take other values of $K$. Note also that on the time interval $(20,30)$ of the most rapid growth of nonpassionaries (constituting the main part of the population) their amount doubles. That corresponds to doubling time  $\frac{10\cdot 900}{60}=150$ years which is roughly   in line with observations on the maximal growth rate of human populations.
\begin{figure}[ht]
	\centering
	\includegraphics[width=1.\textwidth]{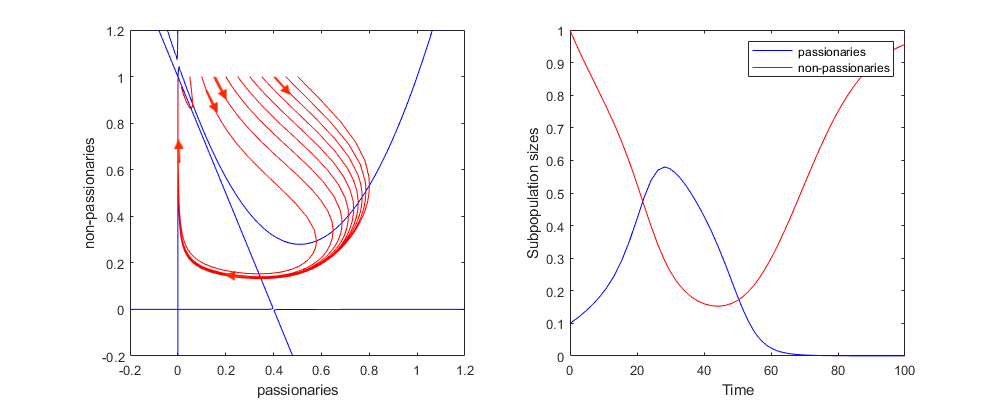}
	\caption{Illustration of the excitable dynamics in the system \eqref{eq.nonlinear1}. Panel A: nullclines $\dot x=0$ and $\dot y=0$ are shown in blue and a set of phase trajectories with starting points  $x_0=0.05,0.1, ..., 0.5$  and $y_0$ are shown in red. Panel B: time evolution of both variables, $x$ (blue line) and $y$ (red line) from the initial condition  $x=0.1$, $y=y_0$. Parameters values: $\alpha=0.02$, $y_0=1$, $\beta_1=1/3$, $\beta_2=-2.5$, $\gamma=0.1$. }
	\label{fig:Null1}
\end{figure} 

For the dynamics illustrated by Figure \ref{fig:Null} it is essential that $\beta_1<0$ and $\beta_2>0$, that is, passionaries are suppressed by the rest of the population while their own impact to non-passionaries is positive. It appears that such relationships between passionaries and non-passionaries  is not the only one allowing excitable dynamics. An alternative case when $\beta_1>0$ and $\beta_2<0$, that is, passionaries are activated by non-passionaries which are in turn suppressed by the passionaries can also result in  the exciatable dynamics. This scenario is illustrated in Figure \ref{fig:Null1} which  similarly to the Figure \ref{fig:Null} shows the dynamics  of  the system \eqref{eq.nonlinear1} but with  swapped signs of the parameters $\beta_1$ and $\beta_2$ and with $y_0=1$ (rather than $y_0=0.05$ in Figure \ref{fig:Null}). Panel A shows null-clines of the system and typical phase trajectories obtained from the over-threshold perturbation of the steady equilibrium   $(x_2,y_2)=(0,y_0)$. Perturbation is made by an increase of the $x$-value over the threshold, $\alpha=0.02$. Time dependence of the variables $x$ and $y$ for one of the phase trajectories is shown on panel B. Now we see that $x$ and $y$ change in the opposite directions: initially $x$ is increasing (for $t<28$) and $y$ is decreasing (for $t<44$) and later  the both variables inverse their rate of change. Furthermore, here $y_0=1$ and the entire dynamics can be seen as change in fraction of passionaries over time for the population of a roughly constant size. Similarly to the case shown in Fig. \ref{fig:Null}, model time unit corresponds to 15 years. 

There is an important difference between the dynamics shown in Figures \ref{fig:Null} and \ref{fig:Null1}. In Figure \ref{fig:Null} the spike in variable $x$ is followed by the spike in the variable $y$, that is, an increase in number of passionaries is followed by the increase of the size of remaining population. As for the dynamics shown in Figure \ref{fig:Null1} we notice that the increase in number of passionaries is followed by the decrease of the size of remaining population. We know that an increase in number of pasionaries results in the expansion of the polity and correspondingly to the growth of the population. Hence the dynamics shown in Figure \ref{fig:Null} looks natural if the variable $y$ represents harmonious people. However if the variable $y$ is associated with subpassionaries then the dynamics in Figure \ref{fig:Null1} is not impossible as subpassionaries may be suppressed by the passionaries. Up to now the variable $y$ was considered as including both, harmonious people and subpassionaries. In order to consider their dynamics separately we will modify our model by allocating variables to each of these two subpopulations.  

\section{Three-Variable Model}

Gumilev in his theory of ethnogenesis considers three types of individuals who constitute ethnos and whose behaviour has an impact to the ethnogenetic process. To follow this concept  we extend the two-variable model described by the system \eqref{eq.nonlinear1} by including extra variable $z$, so that the variables, $x$, $y$ and $z$, represent the sizes of subpopulations of  passionaries ($x$), harmonious people ($y$) and subpassionaries ($z$). Furthermore, we will presume that the dynamics of subpassionaries is similar to that of passionaries that is, their rate of change has quadratic dependence on their own sizes. However passionaries and subpassionaries differ by their relationships with harmonious people and each others. So, our three-variable system can be represented as the following: 
\begin{equation}\label{eq.three-var1}
\left\{
\begin{array}{l}
\dot{x}=\gamma_1 x[(1-x)(x-\alpha_1)+\beta_{12} (y-y_0)+\beta_{13} (z-z_0)]; \\
\dot{y}=\gamma_2 y(y_0 - y+\beta_{21} x+\beta_{23} (z-z_0)); \\
\dot{z}=\gamma_3 z [(z_0-z)(z-\alpha_2) +\beta_{31}x+ \beta_{32} (y-y_0)]
\end{array}
\right.
\end{equation}
with the initial condition $x(0),y(0),z(0)>0$.  
Since its right-hand part is locally Lipschitz, the system (\ref{eq.three-var1}) admits the unique local solution  \cite[Theorem 2.2.]{tes}. 
It looks that in general it can be extended to the unique bounded solution on the infinite horizon $[0,\infty)$.

Excitable dynamics can be observed in the system \eqref{eq.three-var1} under different kinds of interactions between the variables. In general, the steady states of the system (\ref{eq.three-var1}) with $x=0$ are as follows:
\begin{itemize}
\item $(x_1,y_1,z_1)=(0,0,0)$;
\item $(x_2,y_2,z_2)=(0,y_0-\beta_{23}z_0,0)$;
\item $(x_{3,4},y_{3,4},z_{3,4})=(0,y_{3,4},z_{3,4})$, where $y_{3,4}$ and $z_{3,4}$ are the two solutions to equations
$$\left\{\begin{array}{l}
y_0-y+\beta_{23}(z-z_0)=0;\\
(z_0-z)(z-\alpha_2)+\beta_{32}(y-y_0)=0.
\end{array}\right.$$
\item $(x_{5,6},y_{5,6},z_{5,6})=(0,0,z_{5,6})$, where  $z_{5,6}$ are the two solutions to equation
$$(z_0-z)(z-\alpha_2)-\beta_{32}y_0=0.$$
\end{itemize}

\begin{figure}[H]
	\centering
	\includegraphics[width=1.0 \textwidth]{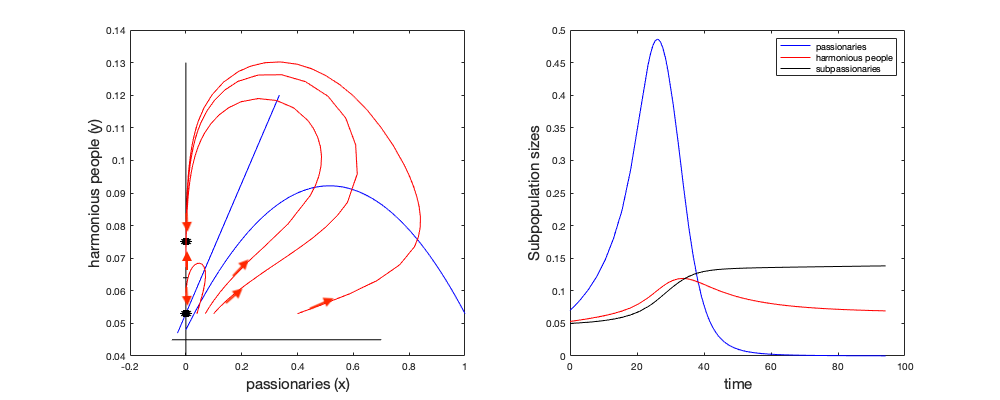}
	\caption{Illustration of the excitable dynamics in the system \eqref{eq.three-var1}. Panel A: nullclines $\dot x=0$ and $\dot y=0$ at $z=z_0$ are shown in blue and a set of phase trajectories with starting points  $x(0)=0.04,0.07, 0.1, 0.4$ and $y(0)=0.053$, $z(0)=0.05$  are shown in red. Panel B: time evolution of all three variables, $x$ (blue line), $y$ (red line) and $z$ (black line) from the initial condition $x(0)=0.07$, $y(0)=0.053$, $z(0)=0.05$. Parameters values: $\alpha_1=0.03$, $\alpha_2=0.11$, $y_0=0.075$, $z_0=0.22$, $\beta_{12}=-6$, $\beta_{13}=0.6$, 
	$\beta_{21}=0.2$, $\beta_{23}=0.1$,
	$\beta_{31}=0.5$, $\beta_{32}=0$,
	$\gamma_1=1$, $\gamma_2=0.7$, $\gamma_3=0.2$.
   	\label{fig:ThreeVar6}}
\end{figure} 

In the case illustrated in Fig. \ref{fig:ThreeVar6}, passionaries are suppressed 
by harmonious people ($\beta_{12}<0$) and promoted by  subpassionaries ($\beta_{13}>0$); harmonious people are promoted by passionaries ($\beta_{21}>0$) as well as by subpassionaries ($\beta_{23}>0$); subpassionaries are promoted by passionaries ($\beta_{31}>0$) and do not depend on harmonious people ($\beta_{32}=0$). For the set of parameters values used in the simulation shown in Fig. \ref{fig:ThreeVar6}, there are eight real-valued steady states, with seven ones among them having non-negative coordinates. All the steady states with non-zero component $x$ are unstable. Among the six steady states enlisted above, $(x_1,y_1,z_1)=(0,0,0)$, $(x_3,y_3,z_3)=(0,0.064,0.11)$, $(x_5,y_5,z_5)=(0,0,0.22)$ and $(x_6,y_6,z_6)=(0,0,0.11)$
are unstable, and the steady states $(x_2,y_2,z_2)=(0,0.053,0)$ and $(x_4,y_4,z_4)=(0,0.075,0.22)$ are stable. In Fig. \ref{fig:ThreeVar6} the stable steady states are shown with big blobs, and the unstable steady state $(x_3,y_3,z_3)=(0,0.064,0.11)$ between them is indicated as the short line.

The excitation appears, starting in the neighbourhood of the  stable point $(x_2,y_2,z_2)$\linebreak$=(0,0.053,0)$ which is shown as the lower blob on the vertical axis. We assigned  the initial values $z(0)=0.05$ (to give a push from the `cemetery' $z=0$), $y(0)=0.053$; $x(0)$ varies from $0.04$ to $0.4$. If the initial push $x(0)$ is below or slightly above the threshold $\alpha_1$, the system quickly returns back to the state $(x_2,y_2,z_2)$. But larger (still small enough) initial perturbation results in the excitation leading to the second stable point $(x_4,y_4,z_4)=(0,0.075,0.22)$ shown as the upper blob on the vertical axis. Depending on the value of $x(0)$, the trajectory approaches the limit $(x_4,y_4,z_4)$ either from below or from above. Even a small over-threshold perturbation $x(0)=0.07>\alpha_1$ grows up to around $x=0.5$ before it relaxes back to $\lim_{t\to\infty} x(t)=0$. Qualitatively, the picture is similar to that presented in Fig. \ref{fig:Null}. Again, one model time unit corresponds to 15 years.
\begin{figure}[ht]
	\centering
	\includegraphics[width=1. \textwidth]{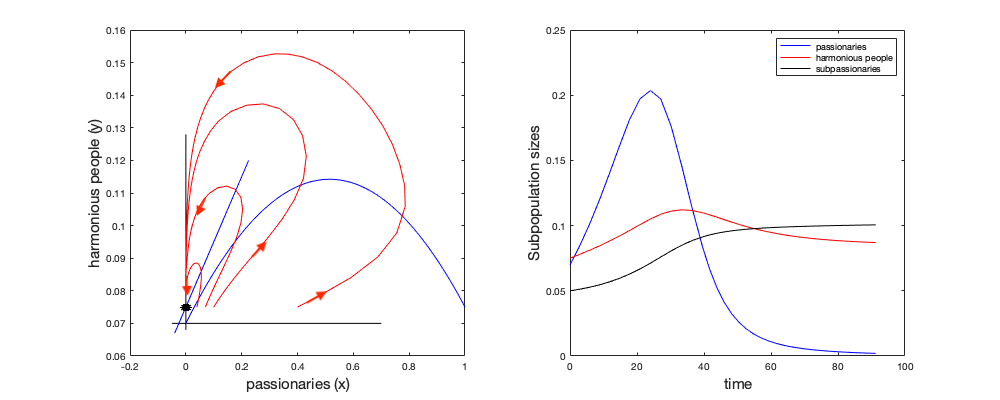}
	\caption{Illustration of the excitable dynamics in the system \eqref{eq.three-var1}. Panel A: nullclines $\dot x=0$ and $\dot y=0$ at $z=z_0$ are shown in blue and a set of phase trajectories with starting points  $x(0)=0.04,0.07,0.1,0.4 $ and $y(0)=0.075$, $z(0)=0.05$  are shown in red. Panel B: time evolution of all three variables, $x$ (blue line), $y$ (red line) and $z$ (black line) from the initial condition $x(0)=0.07$, $y(0)=0.075$, $z(0)=0.05$. Parameters values: $\alpha_1=0.03$, $\alpha_2=0.11$, $y_0=0.075$, $z_0=0.$, $\beta_{12}=-6$, $\beta_{13}=0.6$, 
	$\beta_{21}=0.2$, $\beta_{23}=0.1$,
	$\beta_{31}=0.5$, $\beta_{32}=0$,
	$\gamma_1=1$, $\gamma_2=0.7$, $\gamma_3=0.2$.
   	\label{fig:ThreeVar1}}
\end{figure} 

As a special case, one can put $z_0=0$ leaving the other parameters the same.
For this set of parameters values, used in the simulation shown in Fig. \ref{fig:ThreeVar1}, there is only one stable steady state $(x,y,z)= (0,0.075,0)$: all three previous equilibria $(x_{2,3,4},y_{2,3,4},z_{2,3,4})$ now coincide.
The system again exhibits excitable kinetics: over-threshold perturbation of $x$ ($x>\alpha_1$) grows up to around $x=0.2$ before it relaxes back to the equilibrium.

\begin{figure}[ht]
	\centering
	\includegraphics[width=1.\textwidth]{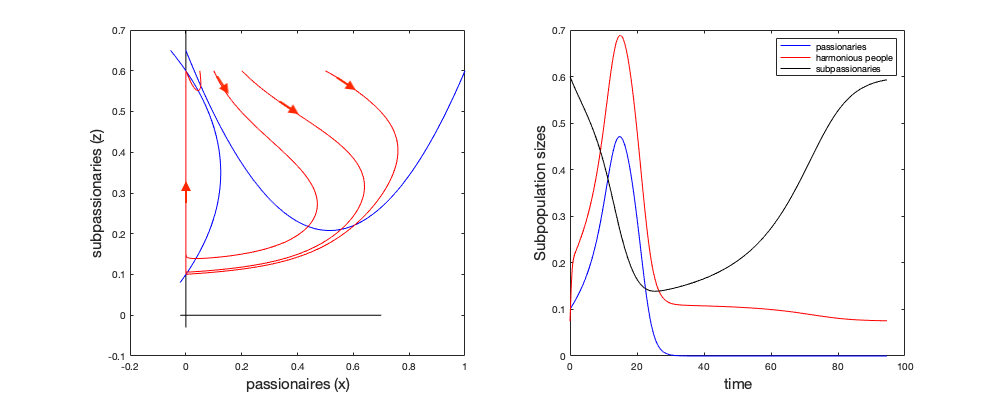}
	\caption{Illustration of the excitable dynamics in the system \eqref{eq.three-var1}. Panel A: nullclines $\dot x=0$ and $\dot z=0$ at $y=y_0$ are shown in blue and a set of phase trajectories with starting points  $x(0)=0.05,0.1,0.2, 0.5$ and $y(0)=y_0,z(0)=z_0$  are shown in red. Panel B: time evolution of all three variables, $x$ (blue line), $y$ (red line) and $z$ (black line) from the initial condition $x(0)=0.1$, $y(0)=y_0$, $z(0)=z_0$. Parameters values: $\alpha_1=0.03$, $\alpha_2=0.1$, $y_0=0.075$, $z_0=0.6$, $\beta_{12}=-0.06$, $\beta_{13}=0.6$, 
		$\beta_{21}=1.25$, $\beta_{23}=-0.075$,
		$\beta_{31}=-0.5$, $\beta_{32}=0$,
		$\gamma_1=2$, $\gamma_2=20$, $\gamma_3=0.6$.
		\label{fig:ThreeVar2}}
\end{figure}

An alternative dynamics for the system \eqref{eq.three-var1} is shown in Fig \ref{fig:ThreeVar2}. Here the interactions between the variables is slightly different from those used for the dynamics illustrated in Fig  \ref{fig:ThreeVar6}. The difference is that harmonious people are  suppressed (rather than promoted) by subpassionaries ($\beta_{23}<0$) and subpassionaries are promoted (rather than suppressed) by passionaries ($\beta_{31}<0$).
For the set of parameters values used in the simulation shown in Fig. \ref{fig:ThreeVar2} there are two stable steady states $(x_2,y_2,z_2)=(0,0.12,0)$ and $(x_4,y_4,z_4)=(0,y_0,z_0)$: the enumeration is in accordance with the expressions below equation (\ref{eq.three-var1}). The system at the  steady state  $(x_4,y_4,z_4)$ exhibits excitable kinetics: over-threshold perturbation of $x$ ($x>\alpha_1$) grows up to around $x=0.7$ before it relaxes back to the equilibrium $(x_4,y_4,z_4)$. As we have two stable states, the relaxation can bring the system to another steady state $(x_2,y_2,z_2)$, and this is observed in the system with slightly different set of parameters values.  Qualitatively, the dynamics illustrated in Fig \ref{fig:ThreeVar2} is similar to that presented in Fig. \ref{fig:Null1}.

One can see that in Fig. \ref{fig:ThreeVar2} the number of subpassionaries dramatically decreases before going back to the equilibrium. This  perhaps does not often occur  in reality. The dynamic in Fig. \ref{fig:ThreeVar1} looks more reasonable. But the growth of the passionaries subpopulation (from $0.07$ to $0.2$) is not as impressive as in Fig. \ref{fig:ThreeVar6} (from $0.07$ to $0.49$). Therefore, in the further simulations, we take the parameters values from the latter case (Fig. \ref{fig:ThreeVar6}). If
 $\beta_{12}\le 0$, $\beta_{32}\le 0$  the existence of the unique bounded solution to the system (\ref{eq.three-var1})
 will follow from the investigation of the stochastic version of the model: see Lemma \ref{l1} and its proof, especially, Remark \ref{rode}.

\section{Ethnogenesis in Noisy Environment}

In order to study the impact of noise to the ethnogenetic process we will modify the three-variable model by adding extra (stochastic) terms to the system (\ref{eq.three-var1}). To justify the modification which we are about to impose, let us consider the following change of variables: 
\begin{equation}\label{e7p}
v_1=\ln x,~~~~~v_2=\ln y,~~~~~v_3=\ln z.
\end{equation}
In these variables the equations (\ref{eq.three-var1}) transform into
\begin{equation}\label{e8}
\left\{
\begin{array}{l}
\dot{v_1}=\gamma_1 [(1-e^{v_1})(e^{v_1}-\alpha_1)+\beta_{12} (e^{v_2}-y_0)+\beta_{13} (e^{v_3}-z_0)]; \\
\dot{v_2}=\gamma_2 [y_0 -e^{v_2} +\beta_{21}  e^{v_1}+\beta_{23} (e^{v_3}-z_0)]; \\
\dot{v_3}=\gamma_3  [(z_0-e^{v_3})(e^{v_3}-\alpha_2) +\beta_{31}e^{v_1}+ \beta_{32} (e^{v_2}-y_0)],
\end{array}
\right.
\end{equation}
with the initial conditions $v_1(0)=\ln x(0)$, $v_2(0)=\ln y(0)$, $v_3(0)=z(0)$. Note that $v_1,v_2$ and $v_3$ may be negative and the initial conditions $x(0),y(0),z(0)>0$ are assumed to be fixed.

The natural way to define the stochastic version is to introduce stochastic differential equations
\begin{equation}\label{e9}
\left\{
\begin{array}{l}
d{V_1}=\gamma_1 [(1-e^{V_1})(e^{V_1}-\alpha_1)+\beta_{12} (e^{V_2}-y_0)+\beta_{13} (e^{V_3}-z_0)]dt+\sigma_1dW_1; \\
d{V_2}=\gamma_2 [y_0 -e^{V_2} +\beta_{21}  e^{V_1}+\beta_{23} (e^{V_3}-z_0)]dt+\sigma_2 dW_2; \\
d{V_3}=\gamma_3  [(z_0-e^{V_3})(e^{V_3}-\alpha_2) +\beta_{31}e^{V_1}+ \beta_{32} (e^{V_2}-y_0)]dt+\sigma_3 dW_3;\\
V_1(0)=v_1(0),~~~~~V_2(0)=v_2(0),~~~~~V_3(0)=v_3(0),
\end{array}
\right.
\end{equation}
which we understand in the sense of the Ito stochastic calculus \cite{oks}. Here $W_1$, $W_2$ and $W_3$ are mutually independent standard Brownian motions on the complete filtered probability space $(\Omega,{\cal F},({\cal F}_t)_{t\ge 0},P)$, and $\sigma_1,\sigma_2,\sigma_3>0$. After that,
$$X=e^{V_1},~~~~~Y=e^{V_2},~~~~~Z=e^{V_3}$$
will be the random processes representing the sizes of the subpopulations of passionaries, harmonious people and subpassionaries respectively. Note  that $X,Y$ and $Z$ satisfy stochastic differential equations
\begin{equation}\label{e30}
\left\{
\begin{array}{l}
dX=\gamma_1 X[(1-X)(X-\alpha_1)+\beta_{12} (Y-y_0)+\beta_{13} (Z-z_0)+\frac{\sigma_1^2}{2}]dt+\sigma_1XdW_1; \\
dY=\gamma_2 Y[y_0 - Y+\beta_{21} X+\beta_{23} (Z-z_0)+\frac{\sigma_2^2}{2}]dt+\sigma_2YdW_2; \\
dZ=\gamma_3 Z [(z_0-Z)(Z-\alpha_2) +\beta_{31}X+ \beta_{32} (Y-y_0)+\frac{\sigma_3^2}{2}]dt+\sigma_3ZdW_3;\\
X(0)=x(0),~~~~~Y(0)=y(0),~~~~~Z(0)=z(0).
\end{array}
\right.
\end{equation}
Derivation of \eqref{e30} can be found in \cite{oks} where it is stated as Theorem 4.2.1.
Here and below, capital letters denote random variables and processes.

In what follows, all the coefficients in  (\ref{eq.three-var1}), (\ref{e8}), (\ref{e9}) and (\ref{e30}) are assumed to be positive apart from $\beta_{12}\le 0$ and $\beta_{32}\le 0$.

\begin{lemma}\label{l1}
Stochastic differential equations (\ref{e9}) (and hence (\ref{e30})) have a  unique strong continuous solution on the time horizon $[0,\infty)$.
\end{lemma}

The proof is presented in the Appendix. It implies that, under positive initial conditions, the ordinary differential equations (\ref{eq.three-var1}) (and hence (\ref{eq.nonlinear1})) have a unique solution such that $x(t),y(t),z(t)>0$ for all $t\ge 0$: see Remark \ref{rode}.

\begin{figure}[ht]
	\centering
	\includegraphics[width=1.\textwidth]{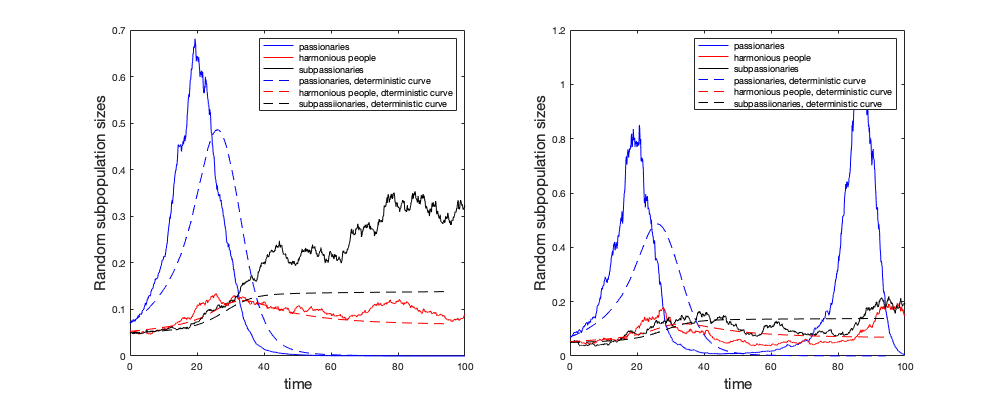}
	\caption{ Stochastic dynamics of the ethnogenesis: examples of the solution to stochastic differential equations (\ref{e30}).
	Parameters values: $\alpha_1=0.03$, $\alpha_2=0.11$, $y_0=0.075$, $z_0=0.22$, $\beta_{12}=-6$, $\beta_{13}=0.6$, 
	$\beta_{21}=0.2$, $\beta_{23}=0.1$,
	$\beta_{31}=0.5$, $\beta_{32}=0$,
	$\gamma_1=1$, $\gamma_2=0.7$, $\gamma_3=0.2$, 	Initial conditions: $X(0)=0.07$, $Y(0)=0.053$, $Z(0)=0.05$. 
Panel A:	
	$\sigma_1=\sigma_2=\sigma_3=0.05$.	
Panel B: $\sigma_1=\sigma_2=\sigma_3=0.1$.
		\label{f101}}
\end{figure} 

Examples of  the stochastic dynamics  exhibited in the system \eqref{e30} are shown in Fig. \ref{f101}. The values of model parameters used for this illustration are identical to those in Fig. \ref{fig:ThreeVar6}. The amplitude of noise in Panel B is twice higher than in Panel A. Solid lines in the both panels show the stochastic dynamics for three subgroups composing the population. For comparison, we also provide the deterministic curves which are represented by the dashed lines (note, that they are identical to the lines shown in Fig. \ref{fig:ThreeVar6}, panel B). It is evident that the stochastic dynamics is significantly different from the deterministic one. One can see from Panel A that the accumulation of noise results in the much higher bust in the level of passionaries: the amplitude of the bust in the stochastic case is about 0.7 against 0.5 in the deterministic one. Also in Panel A we see that the noise causes significant change in the level of subpassionaries: at t=100 this level in stochastic case is over 0.3 while it should be about 0.13 in the deterministic case. The impact of noise is even more evident from Panel B where the amplitude of noise is twice higher than in Panel A. We can see that the increase in the level of noise not only increases the discrepancy with the deterministic case (the amplitude of bust in Panel B is over 0.8) but also results in the occurrence of a new bust. While the first bust in Panel B was initiated manually, the second bust appears due to the stochastic effects in the system. This observation leads us to the conclusion that the ethnogenesis can be initiated by the noise in the environment surrounding the population.

\section{Interaction of Ethnogenetic Processes}

We conclude our study with modelling the interaction between two ethnoses, following the same ethnogenetic processes, which  however are shifted over time. The (random) sizes of subpopulations of passionaries, harmonious people and subpassionaries for the first ethnos are denoted as $X_1$, $Y_1$ and $Z_1$, while for the second ethnos  as $X_2$, $Y_2$ and $Z_2$. We assume that, being isolated, the ethnoses are identical, described by the stochastic differential equations like (\ref{e30}), but influenced by six mutually independent Brownian motions $W_1,W_2,\ldots,W_6$. We introduce the time lag between two ethnogenetic processes, such that the first one starts at $T=0$, the second ethnos appears $T_1$ time units later than the first one, and communication begins $T_2$ time units later, at the time moment $T_1+T_2$. For simplicity, we also assume that communication is only among the passionaries, and they suppress each other.

Therefore, we investigate the following system of six stochastic differential equations

\begin{equation}\label{e31}
\left\{
\begin{array}{rcl}
dX_1&=&\gamma_1 X_1[(1-X_1)(X_1-\alpha_1)+\beta_{12} (Y_1-y_0)+\beta_{13} (Z_1-z_0)+\frac{\sigma_1^2}{2}\\
&&-c_1\II\{t\ge T_1+T_2\}X_2]dt+\sigma_1X_1dW_1; \\
dY_1&=&\gamma_2 Y_1[y_0 - Y_1+\beta_{21} X_1+\beta_{23} (Z_1-z_0)+\frac{\sigma_2^2}{2}]dt+\sigma_2Y_1dW_2; \\
dZ_1&=&\gamma_3 Z_1 [(z_0-Z_1)(Z_1-\alpha_2) +\beta_{31}X_1+ \beta_{32} (Y_1-y_0)+\frac{\sigma_3^2}{2}]dt+\sigma_3Z_1dW_3;\\
X_1(0)&=&x(0),~~~~~Y_1(0)=y(0),~~~~~Z_1(0)=z(0);\\
~~\\
dX_2&=& \II\{ t\ge T_1\}\Bigl\{ \gamma_1 X_2[(1-X_2)(X_2-\alpha_1)+\beta_{12} (Y_2-y_0)+\beta_{13} (Z_2-z_0)+\frac{\sigma_1^2}{2}\\
&&-c_2\II\{t\ge T_1+T_2\}X_1]dt+\sigma_1 X_2 dW_4\Bigr\};\\
dY_2&=&\II\{t\ge T_1\}\Bigl\{\gamma_2 Y_2[y_0 - Y_2+\beta_{21} X_2+\beta_{23} (Z_2-z_0)+\frac{\sigma_2^2}{2}]dt+\sigma_2Y_2dW_5\Bigr\}; \\
dZ_2&=&\II\{t\ge T_1\}\Bigl\{\gamma_3 Z_2 [(z_0-Z_2)(Z_2-\alpha_2) +\beta_{31}X_2+ \beta_{32} (Y_2-y_0)+\frac{\sigma_3^2}{2}]dt+\sigma_3W_6\Bigr\};\\
X_2(0)&=&X_2(T_1)=x(0),~~~~~Y_2(0)=Y_2(T_1)=y(0),~~~~~
Z_2(0)=Z_2(T_1)=z(0).\\
\end{array}
\right.
\end{equation}
The meaning of all the parameters is the same as in the previous models (i.e. model \eqref{e30}). Two new parameters  ($c_1,c_2>0$) define the strength of suppressive interactions between passionaries in the two ethnic groups. This system of stochastic differential equations has a unique strong continuous solution on the time horizon $[0,\infty)$. The proof of this statement is similar to the proof of Lemma \ref{l1}. 

\begin{figure}[ht]
	\centering
	\includegraphics[width=1.\textwidth]{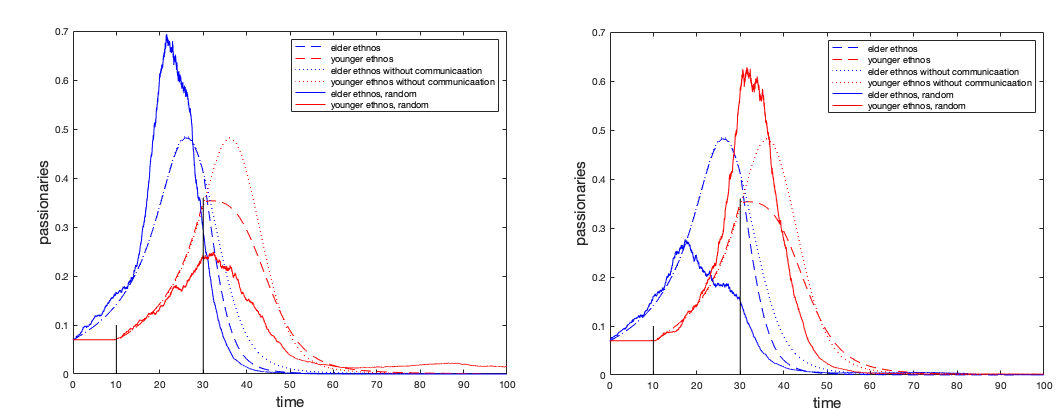}
	\caption{ Stochastic dynamics of comunicating ethnoses.
	Parameters values: $\alpha_1=0.03$, $\alpha_2=0.11$, $y_0=0.075$, $z_0=0.22$, $\beta_{12}=-6$, $\beta_{13}=0.6$, 
	$\beta_{21}=0.2$, $\beta_{23}=0.1$,
	$\beta_{31}=0.5$, $\beta_{32}=0$,
	$\gamma_1=1$, $\gamma_2=0.7$, $\gamma_3=0.2$, $c_1=c_2=0.22$.	 The black vertical lines show the moment of birth of the younger ethnos and the moment when the communication begins. Initial conditions: $X(0)=0.07$, $Y(0)=0.053$, $Z(0)=0.05$.
		\label{f102}}
\end{figure} 
Two examples of dynamics in the interacting ethnoses, described by the system \eqref{e31}, is given in Fig.\ref{f102}. Only the bust dynamics exhibited by the passionaries in  both ethnoses is shown on this figure, with the solid lines showing stochastic dynamics, dashed - deterministic ($\sigma_1=\sigma_2=\sigma_3=0$) and dotted - the deterministic dynamics in the case of non-interacting ethnoses ($c_1=c_2=0$). Dotted lines have identical shapes and this indicates that the two ethnogenetic processes, in the absence of the noise and interaction between the ethnoses, are identical. While the dashed blue line is almost identical to the dotted blue line, the dashed red line is considerably lower than the dotted red line, and this indicates that, in the absence of the noise, the younger ethnos (dashed red line) is suppressed by the older ethnos (dashed blue line). Finally we note that the solid blue line in Panel A is higher than the dashed blue line, while the solid red line is lower than the dashed red line. This observation illustrated the impact of the noise to the dynamics of the interacting ethnoses, which in this particular case results in the amplification of the suppression of the younger ethnos by the older one. Looking at the shapes of  the solid and dashed lines on Panel B we come to the conclusion that the noise can also result in the suppression of the older ethnos by the younger one. Comparing dynamics presented in Panels A and B we note that the dynamics exhibited by two interacting ethnoses is greatly affected by the noise, although, as numerous simulations confirm, the scenario from Panel A is more likely to take place. 

\section{Discussion}

In this work we have presented the mathematical model of ethnogenesis which we have developed on the basis of the paradigm of "passionary tension" introduced by Gumilev \cite{Gumilev1}. According to Gumilev, passionary tension can occur in certain polities as a result of formation and growth of a subgroup of positively motivated people, whom Gumilev called "passionaries". The idea that the growth and evolution of a polity (which can cause the formation of civilisation) is based on its internal structure, and particularly, on the formation of a certain subgroup of people who push the polity forward, was introduced by arabic historian Ibn Khaldun in the 15th  century \cite{khaldun}. One of the main points made by Gumilev is that the measure of passionary tension in the polity is given by its size, i.e. the size of population or territory. As the quantitative data on the territorial expansion and collapse of past civilisations are known much better than on their population sizes, it makes sense to use the size of area taken by a polity as a measure of the passionary tension in this polity.     

The model we have presented here is based on the consideration of the internal structure of the polity with the dynamics of this structure described by ordinary differential equations. The main point about the polity's internal structure is that there is a  subgroup of people, namely, passionaries, and the size of this subgroup gives a measure of the passionary tension in the polity, which in turn can be considered as the measure of the size of territory occupied by the polity. This approach allows to consider the interaction of the given polity with its neighbours indirectly: the polity's geopolitical success is proportional to the number of passionaries in it. 

The main feature of the model we have presented here is that it produces the excitable dynamics in the structure of the evolving polity. That is, when the polity is in equilibrium (in homoeostatic state) there are no passionaries in it. However, if there appear a small number of passionaries, this number grows up to considerable level and then declines back to zero. Thus, formation of busts, describing the raise and fall of civilisations, is considered here as an excitation process. Using different versions of the model we have performed the following studies:
\begin{itemize}
	\item In the two-variable model given by \eqref{eq.nonlinear1} we explored the types of interactions between the passionaries and the rest of the population resulting in the excitable kinetics. Two types of such interactions are illustrated in Figures 2 and 3.
	\item In the three-variable model given by \eqref{eq.three-var1} we explored the types of interactions between three groups allowing to observe the excitable kinetics. Three types of such interactions are illustrated in Figures 4, 5 and 6.
	\item In the stochastic model given by \eqref{e30} we have found that the noise, when imposed into the model equations, tends to amplify and results in significant variations in the amplitude of the bust in the system (see Figure 7). One can conclude that such noise probably adds to the variation of the territorial size and duration of life of different civilisations. 
	\item In the model of interaction polities given by \eqref{e31} we studied the interaction of the polities of different age undergoing ethnogenesis. We found that if they interact in a way that passionaries from one polity suppress the passionaries from the other, then the older polity will be more successful, if the success is measured by the number of passionaries in the polity. However, this is not necessarily the case when we impose the noise (see Figure 8).  
\end{itemize}

The presented model can be extended in various ways for further studies. One of such studies can focus on the interaction of polities under a range of different assumptions about the ways these polities interact. Another obvious direction for future research is to extend the model in order to fit it to available observation data. 

\section{Appendix}

\begin{figure}[ht]
	\centering
	\includegraphics[width=1.\textwidth]{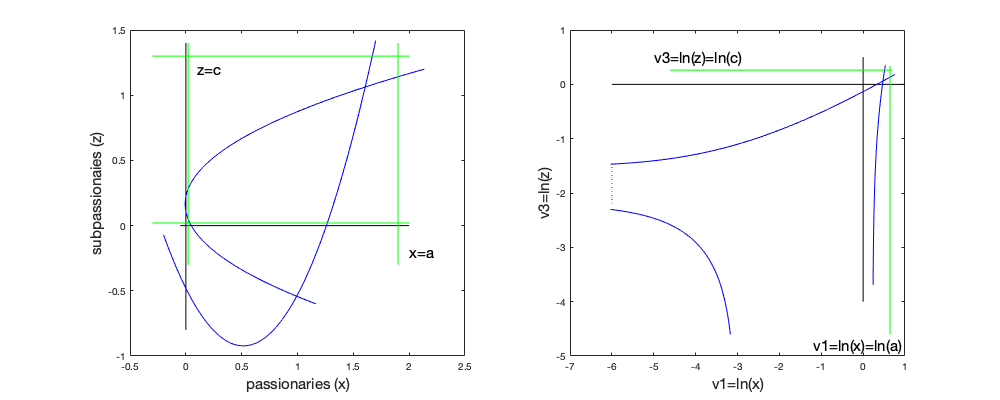}
	
	\caption{Panel A: nullclines $\dot x=0$ and $\dot z=0$ at $y=0$ are shown in blue. Panel B: nullclines $\dot v_1=0$ and $\dot v_3=0$ at $v_2\to-\infty$ are shown in blue. The dotted blue line corresponds to the part of the parabola $\dot z=0$ with negative values of $x$. Parameters values: $\alpha_1=0.03$, $\alpha_2=0.11$, $y_0=0.075$, $z_0=0.22$, $\beta_{12}=-6$, $\beta_{13}=0.6$, 
	$\beta_{21}=0.2$, $\beta_{23}=0.1$,
	$\beta_{31}=0.5$, $\beta_{32}=0$,
	$\gamma_1=1$, $\gamma_2=0.7$, $\gamma_3=0.2$.
		\label{f7}}
\end{figure}

\underline{Proof of Lemma \ref{l1}.} Similarly to  Fig. \ref{fig:ThreeVar6}, panel A, we present the nullclines $\dot x=0$ and $\dot z=0$ at $y=0$ in Fig. \ref{f7}, panel A. Fix a point $(a,c)\in\RR_+\times\RR_+$ such that
\begin{eqnarray*}
(1-a)(a-\alpha_1)-\beta_{12}y_0+\beta_{13}(z-z_0)<0 && \mbox{for all } z\le c;\\
(z_0-c)(c-\alpha_2)+\beta_{31}x-\beta_{32}y_0<0 && \mbox{for all }x\le a.
\end{eqnarray*}
The half-open rectangle $(0,a]\times(0,c]$ is shown with the green lines. Clearly, it is always possible to increase simultaneously $a$ and $c$, so we assume that $a,c>1$, $0<x(0)<a$ and $0<z(0)<c$.
Roughly speaking, the point $(a,c)$ is outside the `internal part' of the both parabolas.

The similar picture in the variables (\ref{e7p}) is given in Fig. \ref{f7}, panel B: the images of the parabolas, shown with the blue lines, represent the nullclines $\dot v_1=0$ and $\dot v_3=0$ of equation (\ref{e8}) in the limiting case when $v_2\to -\infty$.

In the space $\RR^3$, consider the infinite closed prism $\Pi$ defined by 
$$-\infty<v_1\le \ln a,~~~~~-\infty<v_2\le \ln b,~~~~~-\infty<v_3\le \ln c,$$
where 
\begin{equation}\label{e8pr}
b>y_0+\beta_{21}a+\beta_{23}(c-z_0).
\end{equation}
Without loss of generality, we assume that $y(0)<b$ and $b>1$.

\begin{remark}\label{rst1}
For such a prism, we have the following.
\begin{itemize}
\item If $x=a$, then the square bracket in the first equation (\ref{eq.three-var1}) is negative for all $z\le c$ and all $y>0$. (Recall that $\beta_{12}\le 0$.) Therefore, if $v_1= \ln a$, then the square bracket in the first equation (\ref{e8}) and (\ref{e9}) is negative for all $v_3\le\ln c$ and all $v_2>-\infty$.
\item Similarly, if $z=c$ then the square bracket in the third equation (\ref{eq.three-var1}) is negative for all $x\le a$ and all $y>0$. (Recall that $\beta_{32}\le 0$.) Therefore, if $v_3= \ln c$, then the square bracket in the third equation (\ref{e8}) and (\ref{e9}) is negative for all $v_1\le\ln a$ and all $v_2>-\infty$.
\item If $y=b$, then the square bracket in the second equation (\ref{eq.three-var1}) is negative for all $x\le a$ and  $z\le c$. Therefore, if $v_2= \ln b$, then the square bracket in the second equation (\ref{e8}) and (\ref{e9}) is negative for all $v_1\le\ln a$ and all $v_3\le \ln c$.
\end{itemize}
\end{remark}

Now we modify the equations  (\ref{e9})  outside $\Pi$: if $(v_1,v_2,v_3)=m(\hat v_1,\hat v_2,\hat v_3)$ for some $m>1$ with $(\hat v_1,\hat v_2,\hat v_3)\in\partial\Pi$, then we put
$$
\left\{
\begin{array}{l}
f_1(v_1,v_2,v_3):=\gamma_1 [(1-e^{\hat v_1})(e^{\hat v_1}-\alpha_1)+\beta_{12} (e^{\hat v_2}-y_0)+\beta_{13} (e^{\hat v_3}-z_0)]; \\
f_2(v_1,v_2,v_3):=\gamma_2 [y_0 -e^{\hat v_2} +\beta_{21}  e^{\hat v_1}+\beta_{23} (e^{\hat v_3}-z_0)]dt; \\
f_3(v_1,v_2,v_3):=\gamma_3  [(z_0-e^{\hat v_3})(e^{\hat v_3}-\alpha_2) +\beta_{31}e^{\hat v_1}+ \beta_{32} (e^{\hat v_2}-y_0)]dt
\end{array}
\right.
$$
and introduce stochastic differential equations (further, SDEs)
\begin{equation}\label{e10}
\left\{
\begin{array}{l}
d V_1=f_1( V_1, V_2, V_3)dt+\sigma_1 dW_1;\\
d V_2=f_2( V_1, V_2, V_3)dt+\sigma_2 dW_2;\\
d V_3=f_3( V_1, V_2, V_3)dt+\sigma_3 dW_3,\\
 V_1(0)=v_1(0),~~~~~ V_2(0)=v_2(0),~~~~~ V_3(0)=v_3(0).
\end{array}
\right.
\end{equation}
They satisfy all the conditions which guarantee the existence of the unique continuous strong  solution \cite[Remark 14.21]{el} or \cite[Theorem 5.2.1]{oks}: all the functions $f_1,f_2$ and $f_3$ are bounded and Lipschitz in $\RR^3$.

\begin{remark}\label{rode}
If $\sigma_1=\sigma_2=\sigma_3=0$, we have just the system of ordinary differential equations which is uniquely solvable on the time horizon $[0,\infty)$ \cite[Corollary 2.6.]{tes}. This solution $(v_1(t),v_2(t),v_3(t))$ can never leave the prism $\Pi$ because on the bounday $\ln a$ for the component $v_1$ ($\ln b$ for $v_2$ and $\ln c$ for $v_3$) the derivative $\dot v_1$ is negative ($\dot v_2<0$ and $\dot v_3<0$ correspondingly): see Remark \ref{rst1}.

Finally, within the prism $\Pi$, the vector $(v_1(t),v_2(t),v_3(t))$ satisfies differential equations (\ref{e8}) and the functions $x(t):=e^{v_1(t)}$, $y(t):=e^{v_2(t)}$ and $z:=e^{v_3(t)}$ are well defined on the infinite horizon $[0,\infty)$, satisfy equations (\ref{eq.three-var1}) and are strictly positive. As was noted below (\ref{eq.three-var1}), these equations cannot have other solutions.
\end{remark}

In the stochastic version with $\sigma_1,\sigma_2,\sigma_3>0$, the solution to  SDE (\ref{e10})  can exit any one prism on a finite time interval. We need to define a sequence of increasing prisms $\{\Pi_i\}_{i=0}^\infty$ coming from a carefully selected sequence $\{(a_i,b_i,c_i)\}_{i=0}^\infty$. Namely, we require that, for a preliminarily fixed $k>0$, the following condition is satisfied.

\begin{condition}\label{con1}
For each $i>0$ for all
\begin{eqnarray*}
&& a_{i-1}\le x\le a_i,~~~~~0<y\le b_i,~~~~~0<z\le c_i\\
(\mbox{or} && 0<x\le a_i,~~~~~b_{i-1}\le y\le b_i,~~~~~0<z\le c_i\\
\mbox{or} && 0<x\le a_i,~~~~~0< y\le b_i,~~~~~c_{i-1}\le z\le c_i)\\
\end{eqnarray*}
the square bracket in the first (second, third) equation (\ref{eq.three-var1}) is negative and $a_i\ge a_{i-1} e^k$ ($b_i\ge b_{i-1} e^k$, $c_i\ge c_{i-1}e^k$ correspondingly). As the result, for all $(v_1,v_2,v_3)\in\Pi_i$ with $v_1\in[\ln a_{i-1},\ln a_i]$ (with $v_2\in[\ln b_{i-1},\ln b_i]$, $v_3\in[\ln c_{i-1},\ln c_i]$) the square  bracket in the first (correspondingly, second, third) equation (\ref{e8}) and (\ref{e9}) is negative. 

Additionally, $\ln a_i\ge \ln a_{i-1}+k$ ($\ln b_i\ge \ln b_{i-1}+k$, $\ln c_i\ge \ln c_{i-1}+k$) and  $x(0)<a_0,y(0)<b_0,z(0)<c_0$.
\end{condition}

Clearly, under this condition,
$$\inf\{|\vec u_i-\vec u_{i+1}|:~~\vec u_i\in\partial\Pi_i,~\vec u_{i+1}\in\partial\Pi_{i+1}\}\ge k>0.$$

Let us explain why Condition \ref{con1} can be satisfied for an arbitrarily fixed $k>0$.

\begin{figure}[ht]
	\centering
	\includegraphics[width=1.\textwidth]{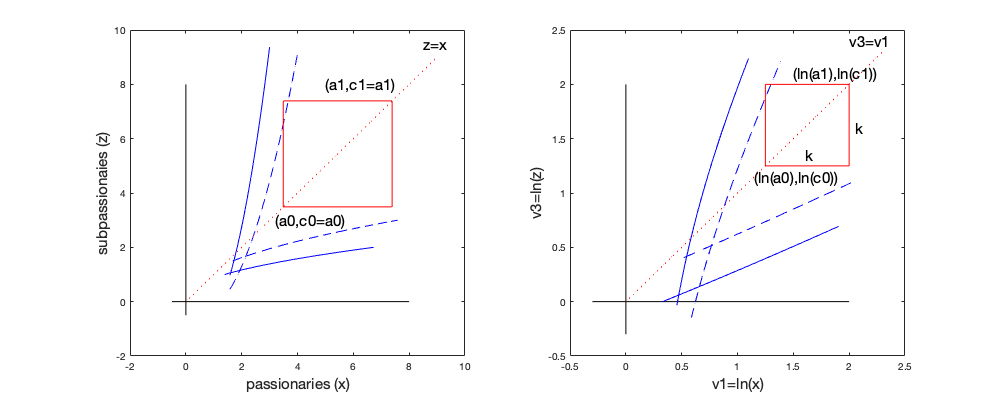}
	\caption{Construction of the first prisms $\Pi_0$ and $\Pi_1$. Parameter values are as in Fig. \ref{f7}; $k=0.75$.
		\label{f8}}
\end{figure}

Along with the parabolas as in Fig. \ref{f7}, panel A, we introduce the expanded graphs (shown in Fig. \ref{f8}, panel A with the dashed blue lines) of the functions
\begin{eqnarray*}
z & = & \frac{-1}{\beta_{13}e^k}\left[(1-x)(x-\alpha_1)-\beta_{12}y_0-\beta_{13}z_0\right];\\
x & = & \frac{-1}{\beta_{31}e^k}\left[(z_0-z)(z-\alpha_2)-\beta_{32}y_0\right].
\end{eqnarray*}
The right-hand parts become bigger than $z$ and $x$ for big enough $x$ and $z$ correspondingly, and one can choose $a_0=c_0$ such that
\begin{eqnarray*}
a_0=c_0 & < &\frac{-1}{\beta_{13}e^k}\left[(1-x)(x-\alpha_1)-\beta_{12}y_0-\beta_{13}z_0\right];\\
a_0=c_0 & < & \frac{-1}{\beta_{31}e^k}\left[(z_0-z)(z-\alpha_2)-\beta_{32}y_0\right]\\
\mbox{and} && a_0>x(0),~~~~~c_0>z(0).
\end{eqnarray*}
After that, the whole red square in Fig. \ref{f8}, panel A, with $a_1=a_0e^k$ and $c_1=c_0e^k=a_1$, is within the  area where the square brackets in the first and third equations (\ref{eq.three-var1}) are negative. The image of Fig. \ref{f8}, panel A on the plain $(v_1,v_3)$ is given in Fig. \ref{f8}, panel B. It remains to take
$$b_0>\max\{y_0+\beta_{21}a_1+\beta_{23}(c_1-z_0),~~~y(0)\}.$$
In general, for all $i\ge 1$, we put 
\begin{eqnarray*}
&& a_i=a_{i-1}e^k,~~~~~c_i=c_{i-1}e^k\\
\mbox{and} && b_i=\max\{y_0+\beta_{21}a_{i+1}+\beta_{23}(c_{i+1}-z_0),~~b_{i-1}e^k\}.
\end{eqnarray*}
The obtained sequence $\{(a_i,b_i,c_i)\}_{i=0}^\infty$ satisfies Condition \ref{con1}.

The SDE (\ref{e10}), for the prism $\Pi_i$, in its vector form means that
\begin{equation}\label{e11}
\vec V^i(t,\omega)=\vec v(0)+\int_0^t\vec f^i(\vec V^i(s,\omega))ds+\Xi\vec W(t,\omega) ~~~~~a.s.
\end{equation}
for all $t\ge 0$. Here the functions $\vec f^i=(f^i_1,f^i_2,f^i_3)$ are constructed, as described above, for the prism $\Pi_i$. The vector notations are conventional, $\Xi=\left[\begin{array}{ccc}
\sigma_1 & 0 & 0\\ 0 & \sigma_2 & 0 \\ 0 & 0 & \sigma_3
\end{array}\right]$. Similarly, we write down equations (\ref{e9}) as
\begin{equation}\label{e12}
\vec V(t,\omega)=\vec v(0)+\int_0^t\vec f(\vec V(s,\omega))ds+\Xi\vec W(t,\omega) ~~~~~a.s.
\end{equation}
for all $t\ge 0$. Note that, if a random process $\vec Z(t,\omega)$ satisfies equation (\ref{e11}), then $\vec V^i(t,\omega)=\vec Z(t,\omega)$ for all $t\ge 0$ almost surely. (The processes $\vec V^i$ and $\vec Z$ are `indistinguishable', often called 'versions' or 'modifications' \cite{oks}.)

Suppose $\tau>0$ is arbitrarily fixed, construct the process $\vec V$ on $[0,\tau]$, which satisfies equation (\ref{e12}), and prove that it is unique. The idea is as follows.{\it 
\begin{itemize}
\item The process $\vec V$ will be a combination of the processes $\vec V^i$.
\item Between the prisms $\Pi_{i-1}$ and $\Pi_i$, the process $\vec V^i$ is pushed back to the prism $\Pi_{i-1}$, and the chance for it to leave the prism $\Pi_{i-1}$ is smaller than a constant $\varepsilon<1$.
\item Hence, almost surely, there is $j\ge 0$ such that the process $\vec V^j$ lives in the prism $\Pi_j$, where it is unique and coincides with $\vec V$.
\end{itemize}}
After that, one can extend the unique solution to (\ref{e12}) to the infinite horizon $[0,\infty)$.

We introduce `debutes'
$$D_i(\omega):=\inf\{0<t\le\tau:~~\vec V^i(t,\omega)\in \Pi_{i+1}\setminus \Pi_i\},~~~i=0,1,\ldots$$
and the (measurable) sets
$$\Omega_i:=\{\omega\in\Omega:~~D_i(\omega)\ge\tau\},~~i=0,1,\ldots.$$
Any debute is a Markov moment because all the processes $\vec V^i$ are continuous a.s. \cite[Example 7.2.2.]{oks}. On each set $\Omega_i$, for all $s\le\tau$, $\vec V^i(s,\omega)\in\Pi_i$ meaning that  $\vec f^i(\vec V^i(s,\omega))=\vec f(\vec V^i(s,\omega))$ and, by (\ref{e11}), for all $t\le\tau$, for $P$-almost all $\omega\in\Omega_i$,
\begin{equation}\label{e20}
\vec V^i(t,\omega)=\vec v(0)+\int_0^t\vec f(\vec V^i(s,\omega))ds+\Xi\vec W(t,\omega).
\end{equation}

For $0\le s\le\tau$, we put
\begin{eqnarray*}
\vec V(s,\omega) &:= & \vec V^0(s,\omega)~\mbox{ on } \Omega_0~~~~~\mbox{ and }\\
\vec V(s,\omega) &:= & \vec V^i(s,\omega)~\mbox{ on }  \overline{\cup_{j=0}^{i-1}\Omega_j}\cap\Omega_i~~~~~\mbox{ for } i=1,2,\ldots.
\end{eqnarray*}
As the result, the continuos process $\vec V$ is built on $\cup_{i=0}^\infty\Omega_i$.

Let us show that 
\begin{equation}\label{e19}
P\left(\bigcup_{i=0}^\infty\Omega_i\right)=1.
\end{equation}
This will be done in three steps.

In what follows, we use expression `on $\Omega'\subset\Omega$ statement$(\omega)$ holds a.s.' to say that $P(\{\omega\in\Omega':~\mbox{statement}(\omega)~\mbox{is false}\})=0$.

\underline{Step 1.} For $i\ge 1$, let us also introduce debutes
$$D_{i,i-1}(\omega):=\inf\{0<t\le\tau:~\vec V^i(t,\omega)\in\Pi_i\setminus\Pi_{i-1}\}$$
and the corresponding subsets
$$\Omega_{i,i-1}:=\{\omega\in\Omega:~~D_{i,i-1}(\omega)\ge\tau\},$$
and show that 
\begin{equation}\label{e13}
P(\Omega_{i-1}\triangle\Omega_{i,i-1})=0.
\end{equation}

Consider the modified $\vec V^{i-1}$ process
$$\vec V^{i-1'}(t,\omega):=\left\{\begin{array}{ll}
\vec V^i(t,\omega), & \mbox{ if } \omega\in\Omega_{i,i-1};\\
\vec V^{i-1}(t,\omega), & \mbox{ if } \omega\in\overline{\Omega_{i,i-1}}
\end{array}\right.$$
and show that 
\begin{equation}\label{e14p}
\vec V^{i-1'}(t,\omega)=\vec V^{i-1}(t,\omega)~~~~~a.s.
\end{equation}
for all $t\in[0,\tau]$.

On the set $\overline{\Omega_{i,i-1}}$, we have, for all $t\in[0,\tau]$,
\begin{equation}\label{e14}
\vec V^{i-1'}(t,\omega)=\vec v(0)+\int_0^t\vec f^{i-1}(\vec V^{i-1'}(s,\omega))ds+\Xi\vec W(t,\omega) ~~~~~a.s.
\end{equation}
because here $\vec V^{i-1'}=\vec V^{i-1}$. On the set $\Omega_{i,i-1}$, $\vec V^{i-1'}=\vec V^{i}\in\Pi_{i-1}$ for all $t\in[0,\tau]$, and, for the vectors $\vec v$ from $\Pi_{i-1}\subset\Pi_i$ we have $\vec f^{i-1}(\vec v)=\vec f^i(\vec v)=\vec f(\vec v)$ meaning that for all $t\in[0,\tau]$ again
\begin{equation}\label{e15}
\vec V^{i-1'}(t,\omega)=\vec v(0)+\int_0^t\vec f^{i-1}(\vec V^{i-1'}(s,\omega))ds+\Xi\vec W(t,\omega) ~~~~~a.s.
\end{equation}
From (\ref{e14}) and (\ref{e15}) we deduce that $\vec V^{i-1'}$ is a solution to the SDE (\ref{e11}) at $i-1$, and assertion (\ref{e14p}) follows. Therefore, on $\Omega_{i,i-1}$, $\vec V^{i-1}=\vec V^i\in\Pi_{i-1}$ for all $t\in[0,\tau]$ almost surely, and $P(\Omega_{i-1}\setminus\Omega_{i,i-1})=0$.

In the similar way, we consider the modified $\vec V^i$ process
$$\vec V^{i'}(t,\omega):=\left\{\begin{array}{ll}
\vec V^{i-1}(t,\omega), & \mbox{ if } \omega\in\Omega_{i-1};\\
\vec V^{i}(t,\omega), & \mbox{ if } \omega\in\overline{\Omega_{i-1}},
\end{array}\right.$$
which, for all $t\in[0,\tau]$, satisfies the SDE (\ref{e11}) at $i$. Therefore, on $\Omega_{i-1}$, $\vec V^{i-1}=\vec V^i\in\Pi_{i-1}$ for all $t\in[0,\tau]$ almost surely, and $P(\Omega_{i,i-1}\setminus\Omega_{i-1})=0$. Equality (\ref{e13}) is proved.

\underline{Step 2.} Suppose $k$ is big enough. In fact, the choice of $k$ depends on $\tau$ and $\max\{\sigma_1,\sigma_2,\sigma_3\}$ only: $k$ must only satisfy inequality
\begin{equation}\label{e16p}
\frac{4}{\sqrt{2\pi}}\int_{d_i}^\infty e^{-\frac{y^2}{2}} dy<\frac{1}{3},~~~~~i=1,2,3,
\end{equation}
where $d_i=\frac{k}{2\sigma_i\sqrt{\tau}}$.

We are going to
show that, for each $i\ge 1$, which is fixed below,
\begin{equation}\label{e15p}
P(\overline{\Omega_i})\le \varepsilon P(\overline{\Omega_{i-1}}),
\end{equation}
where,  $\varepsilon<1$ is some $i$-independent constant.

Clearly, 
$\overline{\Omega_i}\subset \overline{\Omega_{i,i-1}}$ and, for $\omega\in\overline{\Omega_i}$,
$0<D_{i,i-1}(\omega)< D_i(\omega)<\tau$
because $\Pi_{i-1}\subset\Pi_i$ and the process $\vec V^i$ is continuous. We will estimate $P(\overline{\Omega_i}|\overline{\Omega_{i,i-1}})$ assuming that $P(\overline{\Omega_{i,i-1}})>0$. (Otherwise, inequality (\ref{e15p}) is trivial.)

Here and below, usually, all the statements hold $P$-a.s., and all the introduced random variables are defined for $P$-almost all $\omega\in\Omega$, without special remarks.

For $\omega\in\overline{\Omega_{i,i-1}}$, the Markov moment $D_{i,i-1}(\omega)$ is smaller than $\tau$, and we denote $\vec U(\omega):=\vec V^i(D_{i,i-1}(\omega),\omega)\in\partial\Pi_{i-1}$. According to the strong Markov property of $\vec V^i(t,\omega)$ \cite[Theorem 7.2.4.]{oks}, the future behaviour of $\vec V^i(t,\omega)$ on $[D_{i,i-1}(\omega),\tau]$ depends only on $\vec U(\omega)$. Let us estimate $P(\overline{\Omega_i}|\vec U(\omega))$. The set $\overline{\Omega_i}$ is split in three disjoint subsets depending on which component first reaches the boundary $\partial\Pi_i$:
\begin{eqnarray*}
E_1 &=& \{\omega:~V^i_1(D_i(\omega),\omega)=\ln a_i\},\\
E_2 &=& \{\omega:~V^i_2(D_i(\omega),\omega)=\ln b_i\}\setminus E_1,\\
\mbox{and}~~E_3 &=& \{\omega:~V^i_3(D_i(\omega),\omega)=\ln c_i\}\setminus(E_1\cup E_2).
\end{eqnarray*}
We shall prove that, for some $i$-independent constant $\delta<\frac{1}{3}$,
\begin{equation}\label{e16}
P(E_j|U(\omega))\le\delta,~~~~~j=1,2,3.
\end{equation}

Suppose $j=1$:
the reasoning for $j=2$ and $j=3$ is similar.

\begin{figure}[ht]
	\centering
	\includegraphics[width=1.\textwidth]{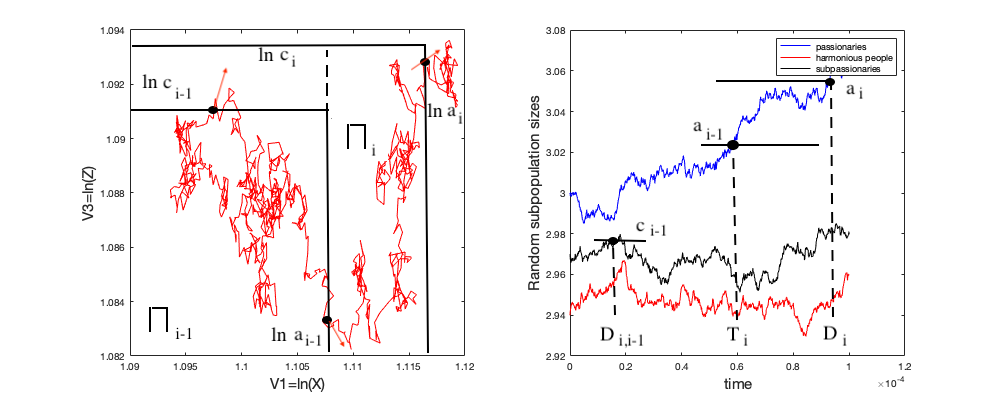}
	\caption{Small fragment of the motion of $\vec V^i$: phase trajectories on panel  A and time trajectories on Panel  B. The black blobs correspond to the time moments $D_{i,i-1}$, $T$ and $D_i$. The red arrows at the blobs on panel A  show the direction of the movement of the phase trajectories at the time moments $D_{i,i-1}$, $T_i$ and $D_i$. 
	Parameters values: $\alpha_1=0.03$, $\alpha_2=0.11$, $y_0=0.075$, $z_0=0.22$, $\beta_{12}=-6$, $\beta_{13}=0.6$, 
	$\beta_{21}=0.2$, $\beta_{23}=0.1$,
	$\beta_{31}=0.5$, $\beta_{32}=0$,
	$\gamma_1=1$, $\gamma_2=0.7$, $\gamma_3=0.2$, $\sigma_1=\sigma_2=\sigma_3=1$.
		\label{f9}}
\end{figure}

Let $\vec u\in\partial\Pi_{i-1}$ be fixed and consider the process $\vec V^i(t,\omega)$  with $\omega\in\overline{\Omega_{i,i-1}}$, $t\ge D_{i,i-1}(\omega)$ as starting from $\vec V^i(D_{i,i-1}(\omega))=\vec u$:
\begin{equation}\label{e17}
\vec V^i(t,\omega)=\vec u+\int_{D_{i,i-1}(\omega)}^t \vec f^i(\vec V^i(s,\omega))ds+\Xi[\vec W(t,\omega)-\vec W(D_{i,i-1}(\omega),\omega)],~t\in [D_{i,i-1}(\omega),\tau].
\end{equation}
Below, we assume that $\omega\in E_1$ in order to estimate $P(E_1|\vec u)$, so that $D_i(\omega)<\tau$ is well defined. The set
$$H_i(\omega):=\{t\in[D_{i,i-1}(\omega),D_i(\omega)]:~V^i_1(t,\omega)\le \ln a_{i-1}\}$$
is not empty: $V^i_1(D_{i,i-1}(\omega))=u_1\le\ln a_{i-1}$. We put 
$$T_i(\omega):=\sup\{t:~t\in H_i(\omega)\}.$$
Although $T_i$ is not a Markov moment, it is a measurable random variable: for each $x\in\RR_+$
\begin{eqnarray*}
\{\omega:~T_i(\omega)>x\}&=&\{\omega:~x< D_i(\omega) \mbox{ and } \exists s\in(x,D_i(\omega))\cap\QQ:~V^i_1(s,\omega)\le a_{i-1}\}\\
&=& \{\omega:~D_i(\omega)> x\}\cap \bigcup_{s\in\QQ,~s>x}\{\omega:~D_i(\omega)>s \mbox{ and }V^i_1(s,\omega)\le \ln a_{i-1}\},
\end{eqnarray*}
where $\QQ$ is the set of rational numbers. Note that $V^i_1(T_i(\omega),\omega)=\ln a_{i-1}$ and remember that $T_i(\omega)\in[D_{i,i-1}(\omega),D_i(\omega)]$ and $0<D_{i,i-1}(\omega)<D_i(\omega)<\tau$. On the time interval $[T_i(\omega),D_i(\omega)]$, the process $\vec V^i(t,\omega)$ is still in $\Pi_i$ and $V^i_1(t,\omega)\in[\ln a_{i-1},\ln a_i]$: see Fig. \ref{f9}.

According to (\ref{e17}) and taking into account Condition \ref{con1}, we have
\begin{eqnarray*}
V^i_1(D_i(\omega),\omega) &=& u_1+\int_{D_{i,i-1}(\omega)}^{T_i(\omega)} f^i_1(\vec V^i(s,\omega))ds+\sigma_1[W_1(T_i(\omega),\omega)-W_1(D_{i,i-1}(\omega),\omega)]\\
&&+\int_{T_i(\omega)}^{D_{i}(\omega)} f^i_1(\vec V^i(s,\omega))ds+\sigma_1[W_1(D_{i}(\omega),\omega)-W_1(T_i(\omega),\omega)]\\
&=& V^i_1(T_i(\omega),\omega)+\int_{T_i(\omega)}^{D_{i}(\omega)} f^i_1(\vec V^i(s,\omega))ds+\sigma_1[W_1(D_{i}(\omega),\omega)-W_1(T_i(\omega),\omega)]\\
&\le & V^i_1(T_i(\omega),\omega)+\sigma_1[W_1(D_{i}(\omega),\omega)-W_1(T_i(\omega),\omega)]:
\end{eqnarray*}
the function $f^i_1(v)$ is negative for $v\in\Pi_i$ with $v_1\in[\ln a_{i-1},\ln a_i]$. Since $V^i_1(T_i(\omega),\omega)=\ln a_{i-1}$ and $V^i_1(D_i(\omega),\omega)=\ln a_i$, we see that, for $\omega\in E_1$,
$$W_1(D_i(\omega),\omega)-W_1(T_i(\omega),\omega)\ge\frac{\ln a_i-\ln a_{i-1}}{\sigma_1}\ge \frac{k}{\sigma_1}.$$
That means 
$$\sup_{0\le s\le\tau} W_1(s,\omega)-\inf_{0\le s\le\tau} W_1(s,\omega)\ge\frac{k}{\sigma_1},$$
and $P(\tilde\Omega)$, the probability of the set $\tilde\Omega$ of all the points $\omega$ satisfying this property, is smaller than some $0<\delta<1/3$ for big enough $k$. To be more specific,
$$\tilde\Omega\subset\{\omega:~\sup_{0\le s\le\tau} W_1(s,\omega)\ge \frac{k}{2\sigma_1} \mbox{ or } \inf_{0\le s\le\tau} W_1(s,\omega)\le \frac{-k}{2\sigma_1}\},$$
and the probability of the set on the right is well studied: it is smaller than
$$4P(\{\omega:~W_1(1,\omega)\ge d\})=\frac{4}{\sqrt{2\pi}}\int_d^\infty e^{-\frac{y^2}{2}}dy,$$
where $d=\frac{k}{2\sigma_1\sqrt{\tau}}$, and approaches zero as $k\to\infty$. See also Lemma \ref{la2}. Thus, for fixed $\vec u$, if $k$ is big enough (and certainly dependent on $\max\{\sigma_1,\sigma_2,\sigma_3\}$ and $\tau$ only), then $P(E_1|\vec u)=\delta<\frac{1}{3}$: inequality (\ref{e16}) is proved.

Applying the similar reasoning to $E_2$ and $E_3$, we conclude that $P(\overline{\Omega_i}|\vec u)\le\varepsilon:=3\delta<1$ and
$$P(\overline{\Omega_i})=\int_{\overline{\Omega_{i,i-1}}} P(\overline{\Omega_i}|\vec U(\omega))dP(\omega)\le \varepsilon P(\overline{\Omega_{i,i-1}}).$$
Inequality (\ref{e15p}) now follows from (\ref{e13}).

\underline{Step 3.} Using (\ref{e15p}), we have
$$P\left(\bigcap_{i=0}^\infty \overline{\Omega_i}\right)=\lim_{N\to\infty}P\left(\bigcap_{i=0}^N \overline{\Omega_i}\right)\le\lim_{N\to\infty} P(\overline{\Omega_N})\le \lim_{N\to\infty}\varepsilon^N P(\overline{\Omega_0})=0.$$
Therefore, $P\left(\bigcup_{i=0}^\infty\Omega_i\right)=1-P\left(\bigcap_{i=0}^\infty \overline{\Omega_i}\right)=1$, and equality (\ref{e19}) is proved.

According to (\ref{e20}), the constructed continuous process $\vec V(t,\omega)$ satisfies equation (\ref{e12}) (and (\ref{e9})) for all $t\in[0,\tau]$ a.s., i.e., it is a strong continuous solution to those SDEs.

Now let us show that, if $\vec Z(t,\omega)$ is a strong continuous solution to the SDE (\ref{e12}), then, for all $t\in[0,\tau]$, $\vec Z(t,\omega)=\vec V(t,\omega)$ a.s.

Since the process $\vec Z$ is continuous a.s., it is bounded a.s. and hence (almost surely) there is (a unique) integer $I(\omega)\ge 0$ such that, for all $t\in[0,\tau]$, $\vec Z(t,\omega)\in\Pi_{I(\omega)}$ and, for $I(\omega)>0$, there is $t\in[0,\tau]$ such that $\vec Z(t,\omega)\in\Pi_{I(\omega)}\setminus \Pi_{I(\omega)-1}$. Recall that $\bigcup_{i=0}^\infty\Pi_i=\RR^3$, so that,  the sets $\tilde\Omega_i:=\{\omega:~I(\omega)=i\}$ are such that $P\left(\Omega\setminus \bigcup_{i=0}^\infty \tilde\Omega_i\right)=0$. On the set $\tilde\Omega_i$, the process $\vec Z$ satisfies SDE
$$\vec Z(t,\omega)=\vec v(0)+\int_0^t \vec f^i(\vec Z(s,\omega))ds+\Xi\vec W(t,\omega),~~~~~t\in[0,\tau]$$
for $P$-almost all $\omega\in\tilde\Omega_i$ because, within the prism $\Pi_i$, $\vec f(\vec v)=\vec f^i(\vec v)$. We see that $\vec Z(t,\omega)=\vec V^i(t,\omega)$ for all $t\in[0,\tau]$ a.s. on $\tilde\Omega_i$ and, for $P$-almost all $\omega\in\tilde\Omega_i$, $\vec V^i(t,\omega)\in\Pi_i$ for all $t\in[0,\tau]$ and, in case  $i>0$, there is $t\in[0,\tau]$ such that $\vec V^i(t,\omega)\in\Pi_i\setminus\Pi_{i-1}$. The last assertion means that, for $P$-almost all $\omega\in\tilde\Omega_i$, $\omega\in\Omega_i$ and, in case $i>0$, $\omega\in\overline{\Omega_{i,i-1}}$. (The process $\vec V^i$ left  $\Pi_{i-1}$ at some moment.) Since $\overline{\Omega_{i,i-1}}\triangle\overline{\Omega_{i-1}}=\Omega_{i,i-1}\triangle\Omega_{i-1}$, for $i>0$, according to (\ref{e13}), $\omega\in\overline{\Omega_{i-1}}$ for $P$-almost all $\omega\in\tilde\Omega_i$. We conclude that, for $P$-almost all $\omega\in\tilde\Omega_i$, $\omega\in\overline{\Omega_{i-1}}\cap\Omega_i$ ($\omega\in\Omega_0$ in case $i=0$) and $\vec V^i(t,\omega)=\vec V(t,\omega)$ for all $t\in[0,\tau]$ a.s. on $\tilde\Omega_i$. As a result, $\vec Z(t,\omega)=\vec V(t,\omega)$ for all $t\in[0,\tau]$ a.s on $\tilde\Omega_i$. Therefore, $\vec Z(t,\omega)=\vec V(t,\omega)$ for all $t\in[0,\tau]$ a.s. on $\bigcup_{i=0}^\infty\tilde\Omega_i$, and the latter set, as explained above, coincides with $\Omega$ up to a set of $P$-measure zero.

The uniqueness of the strong continuous solution to SDE (\ref{e12}) (hence, (\ref{e9})) is proved.

Finally, extension to the infinite horizon $[0,\infty)$ of the solution $\vec V$ to the SDE (\ref{e12}) (and (\ref{e9})) is trivial. Take an increasing sequence $\{\tau_j\}_{j=1}^\infty$, $\tau_j>0$, with $\lim_{j\to\infty}\tau_j=\infty$, construct the solutions $\vec V^{\tau_j}$ to SDE (\ref{e12}) (and (\ref{e9})) on the intervals $[0,\tau_j]$ and put
$$\vec V(t,\omega):=\vec V^{\tau_j}(t,\omega)\II\{t\in[\tau_{j-1},\tau_j)\},$$
where $\tau_0:=0$. Due to the uniqueness of each process $\vec V^{\tau_j}$, for each $k=1,2,\ldots$, $\vec V(t,\omega)=\vec V^{\tau_k}(t,\omega)$ for all $t\in[0,\tau_k)$ a.s. and hence the process $\vec V$ satisfies SDE (\ref{e12}) (and (\ref{e9})) for all $t\in[0,\tau_k)$ a.s. meaning that $\vec V$ satisfies SDE (\ref{e12}) (and (\ref{e9})) for all $t\in[0,\infty)$ a.s.

If there is another process $\vec Z$ satisfying this property then, again due to the uniqueness of $\vec V^{\tau_j}$, $\vec Z(t,\omega)=\vec V^{\tau_j}(t,\omega)$ for all $t\in[\tau_{j-1},\tau_j)$ a.s., $j=1,2,\ldots$. Hence  $\vec Z(t,\omega)=\vec V(t,\omega)$ for all $t\in[0,\infty)$ a.s.

The proof is completed.

\begin{lemma}\label{la2}
Let $W(t,\omega)$ be a Brownian motion on a filtered probability space $(\Omega,{\cal F},({\cal F}_t)_{t\ge 0},P)$. Then, for $a,\tau>0$
$$P(\sup_{0\le t\le\tau}W(t,\omega)\ge a \mbox{ or } \inf_{0\le t\le\tau}W(t,\omega)\le -a)\le\frac{4}{\sqrt{2\pi}}\int_{a/\sqrt{\tau}}^\infty e^{-\frac{y^2}{2}} dy.$$ 
\end{lemma}

\underline{Proof.} Let $T_a(\omega):=\inf\{t:~W(t,\omega)=a\}$. Then, by the symmetry,
\begin{equation}\label{ea1}
P(\sup_{0\le t\le\tau}W(t,\omega)\ge a \mbox{ or } \inf_{0\le t\le\tau}W(t,\omega)\le -a)\le P(T_a(\omega)\le\tau)+P(T_{-a}(\omega)\le\tau)=2P(T_a(\omega)\le\tau).
\end{equation}
Now, by the reflection principle,
$P(T_a\le\tau)=2P(W(\tau,\omega)\ge a)$, and we continue (\ref{ea1}):
$$2P(T_a(\omega)\le\tau)=4\int_a^\infty \frac{1}{\sqrt{2\pi\tau}} e^{-\frac{z^2}{2\tau}}dz=\frac{4}{\sqrt{2\pi}}\int_{a/\sqrt{\tau}}^\infty e^{-\frac{y^2}{2}} dy.$$ 

The proof is completed.


\end{document}